\documentclass{article}
\usepackage{xy}
\usepackage{epsfig}
\usepackage[T1]{fontenc} 
\usepackage[latin1]{inputenc}
\usepackage{amsfonts,amssymb,amsmath}
\usepackage{hyperref} 
\usepackage{graphics} 
\usepackage{multido, xcolor} 

\newtheorem{theo}{\bf\large Theorem}[section]

\newtheorem{rema}[theo]{\bf Remark}

\newtheorem{conj}[theo]{\bf Conjecture}

\def\qed{ \ \hfil$\square$}

\makeindex

\newcommand{\ord}{{\rm ord}}
\newcommand{\SL}{{\rm SL}}

\newcommand{\N}{{\mathbb N}}
\newcommand{\Q}{{\mathbb Q}}

\newcommand{\Z}{{\mathbb Z}}

\newcommand{\e}{\varepsilon}

\font\teneusm=eusm10 \font\seveneusm=eusm7 
\font\fiveeusm=eusm5 
\newfam\eusmfam 
\textfont\eusmfam=\teneusm 
\scriptfont\eusmfam=\seveneusm 
\scriptscriptfont\eusmfam=\fiveeusm

\def\mat #1,#2,#3,#4,{\left({#1\atop #3}{#2\atop #4}\right)}
\def\bra#1,{{\left\lbrace {#1}\right\rbrace}}

\def \Sp{{\rm Sp}}

\def\l1{\langle}

\newcommand{\B}{\left(\begin{array}{cc}}
\newcommand{\E}{\end{array}\right)}

\usepackage{euscript}
\let\scr=\EuScript

\let\mathcal=\scr           
\def\ang#1,{{\left\langle {#1}\right\rangle}}

\def\qed{ \ \hfil$\square$}

\font\teneusm=eusm10 \font\seveneusm=eusm7 
\font\fiveeusm=eusm5 
\newfam\eusmfam 
\textfont\eusmfam=\teneusm 
\scriptfont\eusmfam=\seveneusm 
\scriptscriptfont\eusmfam=\fiveeusm

\font\tengothic=eufm10
\font\sevengothic=eufm7
\font\fivegothic=eufm5
\newfam\Gothic
\textfont\Gothic=\tengothic
\scriptfont\Gothic=\sevengothic
\scriptscriptfont\Gothic=\fivegothic

\def\mat #1,#2,#3,#4,{\left({#1\atop #3}{#2\atop #4}\right)}
\def\bra#1,{{\left\lbrace {#1}\right\rbrace}}

\def \Sp{{\rm Sp}}

\def\l1{\langle}

\def\Numero{${\rm N\sp\circ}$}

\let\scr=\EuScript

\let\mathcal=\scr           

\def\vin{{ {\tiny \mid }  
\kern-7.29pt 
\bigcup }}

\def\ang#1,{{\left\langle {#1}\right\rangle}}

\normalsize

\newcommand{\CC}{\mathbb C}

\newcounter{ncours}{\setcounter{ncours} {1}}

\def\SL {\mathop{\rm SL}\nolimits}
\def\Sp {\mathop{\rm Sp}\nolimits}

\font\tenrus=wncyr10
\font\sevenrus=wncyr7
\newfam\Rus
\textfont\Rus=\tenrus
\scriptfont\Rus=\sevenrus

\def\scr{\fam\Rus\sevenrus} 

\font\tenbrus=wncyb10
\font\eightbrus=wncyb8
\newfam\Brus
\textfont\Brus=\tenbrus
\scriptfont\Brus=\eightbrus

\font\tenirus=wncyi10
\font\eightirus=wncyi8
\newfam\Irus
\textfont\Irus=\tenirus
\scriptfont\Irus=\eightirus

\font\tenrus=wncyr10
\font\sevenrus=wncyr7
\newfam\Rus
\textfont\Rus=\tenrus
\scriptfont\Rus=\sevenrus

\def\scr{\fam\Rus\sevenrus} 
\input cyracc.def
\def\i1{\accent'044i}
\def\I1{\accent'044I}
\def\e1{\accent'040e}
\def\l1{l{}p1}
\newcommand{\T}{{\mathbf{T}}}

\newcommand{\Ta}{{\mathbf{T}(p)}}
\newcommand{\Tb}{{\mathbf{T}_1(p^2)}}
\newcommand{\Tc}{{\mathbf{T}_2(p^2)}}

\newcommand{\Te}{{[\mathbf{p}]}}

\title{
On Zeta Functions and Families of Siegel Modular Forms 
}

\author 
{Alexei PANCHISHKIN 
}
\date{}

\begin{document}
\maketitle


%
%

\begin{abstract}
Let   $p$ be a prime, and let $\Gamma=\Sp_g(\Z)$ be the   Siegel modular group of genus  $g$.
We study  $p$-adic families of  zeta functions and Siegel modular forms. 


$L$-functions of Siegel modular forms are described in terms of 
motivic $L$-functions attached to 
$\Sp_g$, and their analytic properties are given.
 Critical values  for the spinor $L$-functions and $p$-adic constructions are discussed.
Rankin's lemma of higher genus is established.
A general conjecture on a lifting from $\,GSp_{\,2m} \times GSp_{\,2m}$ to $GSp_{\,4m}$ (of genus $g=4m$) is formulated.

Constructions of  $p$-adic families of Siegel modular forms
are given using Ikeda-Miyawaki constructions.


\end{abstract}
  \tableofcontents
\section{Introduction}

The present artical is based on some recent talks of the author, especially the talk in Moscow University
 on February  2, 2007, 
for the International conference 
"Diophantine and Analytical Problems in Number Theory" 
for  100th birthday of  A.O. Gelfond. 

 Let  $\Gamma=\Sp_g(\Z)\subset \SL_{2g}(\Z)$ be the   Siegel modular group of genus  $g$. 
Let   $p$ be a prime,   $\mathbf{T}(p)=\mathbf{T}(\underbrace{1, \cdots, 1}_{g}, \underbrace{p, \cdots, p}_{g})$
the Hecke $p$-operator, and 
 $[\mathbf{p}]_g=p\mathbf{I}_{2g}=\mathbf{T}(\underbrace{p, \cdots, p}_{2g})$ the scalar Hecke operator for $Sp_g$.

We discuss the following topics:

\begin{enumerate}
\item[1)] 
$L$-functions of Siegel modular forms 
\item[2)] 
Motivic $L$-functions for
$\Sp_g$, and their analytic properties
\item[3)] 
 Critical values and $p$-adic constructions for the spinor $L$-functions
\item[4)]{Rankin's Lemma of higher genus}
\item[5)] 
A lifting from $\,GSp_{\,2m} \times GSp_{\,2m}$ to $GSp_{\,4m}$ (of genus $g=4m$)
\item[6)]
Constructions of  $p$-adic families of Siegel modular forms
\item[7)]
Ikeda-Miyawaki constructions and their $p$-adic versions
\end{enumerate}

\section{$L$-functions of Siegel modular forms}

\subsection*{The Fourier expansion of a Siegel modular form.
}
Let  $\displaystyle f=\sum_{\mathcal{T}\in B_n} a(\mathcal{T})q^\mathcal{T}\in \mathcal{M}_k^n$ be a Siegel modular form of weight $k$ and of genus $n$
on the Siegel upper-half plane $\mathbb{H}_n=\{z\in \mathrm{M}_n(\mathbb{C}) \ |\ \mathrm{Im}(z) >0 \}$.

\

Recal some basic facts about 
the formal  {\it  Fourier expansion }  of $f$, which  uses the symbol
$$
q^{\mathcal T}=\exp(2\pi i \, {\rm tr}({\mathcal T}z)) =
 \prod_{i=1}^nq_{ii}^{{\mathcal T}_{ii}}\prod_{i<j}q_{ij}^{2{\mathcal T}_{ij}}
$$
$ \in \mathbb{C}[\![q_{11}, \dots , q_{nn}]\!][q_{ij}, \ q_{ij}^{-1}]_{i,j=1, \cdots, m}$, where
$q_{ij}= \exp(2\pi (\sqrt {-1} z_{i,j}))$,
and
$\mathcal{T}$ is in the   semi-group
 $
B_n =\{{\mathcal T}={}^t{\mathcal T}\ge 0\vert {\mathcal T}\hbox{ half-integral}\}.
$
\subsection*{Hecke operators and the spherical map}
Recall that the local Hecke algebra over $\Z$, denoted by  $\mathcal{L}_{n,\Z}=\Z[\Ta,\Tb,\dotsc,\T_n(p^2)]$, 
 is generated by the following $n+1$ Hecke operators
\begin{align*}
&\Ta:=T(\underbrace{1,\dotsc,1}_{n},\underbrace{p,\dots,p}_{n})\,,\\
&\T_i(p^2):=T(\underbrace{1,\dotsc,1}_{n-i},\underbrace{p,\dotsc,p}_{i},\underbrace{p^2,\dotsc,p^2}_{n-i},\underbrace{p,\dotsc,p}_{i})\,,\,i=1,\dotsc,n
\end{align*}
$\T_n(p^2)=\Te=\Te_n=T(\underbrace{p,\dotsc,p}_{2n})=p\/\mathbf{I}_{2n}$
and the spherical map 
$\Omega:\mathcal{L}_{n,\Q}=\Q[\Ta,\Tb,\dotsc,\T_n(p^2)]\to\Q[x_0,x_1,\dotsc,x_n]$
is an a certain injective ring morphism.
 
\subsection*{Satake parameters of an eigenfunction of Hecke operators}
Suppose that   $f \in{\mathcal M}_k^n$ an eigenfunction of all  {\it Hecke operators}
$f \longmapsto f|T$, $T\in {\mathcal L}_{n,p}$ for all primes $p$,
hence $f|T = \lambda_f(T) f$.

Then all the numbers
$\lambda_f(T) \in \mathbb{C}$ define a homomorphism
$\lambda_f : {\mathcal L}_{n,p} \longrightarrow \mathbb{C}$ 
given by a  $(n + 1)$-tuple of complex numbers
$
    (\alpha_0, \alpha_1, \cdots, \alpha_n) =
     (\mathbb{C}^{\times})^{n+1}  
$
(the Satake parameters of $f$), 
in such a way that 
$$
\lambda_f(T)= \Omega(T)(\alpha_0, \alpha_1, \cdots, \alpha_n).
$$
In particular, one has 
\begin{align*}& 
 \lambda_f(\Te)  =  \alpha_0^2 \alpha_1 \cdots \alpha_n=  p^{kn-n(n+1)/2}  \\ &
 \lambda_f(\Ta)  = \Omega(\Ta)=\alpha_0(1+\alpha_1)\dotsb(1+\alpha_n)=\sum\limits_{j=0}^n\alpha_0s_j(\alpha_1,\alpha_2,\dotsc,
\alpha_n)\,.
\end{align*}


\subsection*
{$L$-functions, functional equation and motives for   $\mathrm{Sp}_n$ (see \cite{Pa94}, \cite{Y})
}
One defines
\begin{itemize}
\item
$\displaystyle
    Q_{f,p}(X)  =  (1 - \alpha_0 X) \prod_{r=1}^n
    \prod_{1 \leq i_1 < \cdots < i_r \leq n} (1 - \alpha_0 \alpha_{i_1}
    \cdots \alpha_{i_r} X),  $

\item
 $\displaystyle R_{f,p}(X)  =  (1 -  X)\prod_{i=1}^n (1 - \alpha_i^{-1} X)
    (1 - \alpha_i X) \in \mathbb{Q}[ \alpha_0^{\pm 1}, \cdots, \alpha_n^{\pm 1}][X].$
\end{itemize}

Then the spinor  $L$-function $L(Sp(f), s)$ and the standard $L$-function $L(St(f), s, \chi)$ of $f$ (for $s\in \mathbb{C}$, and for
all Dirichlet characters) $\chi$ are defined as the Euler products  
\begin{itemize}
\item
$ \displaystyle L(Sp(f), s, \chi) = \prod_{p}Q_{f,p}(\chi(p)p^{-s})^{-1} 
$
\item
$\displaystyle
 L(St(f), s, \chi)    =    \prod_{p}R_{f,p}(\chi(p)p^{-s})^{-1}
$
\end{itemize}

\section{Motivic $L$-functions for
$\Sp_n$, and their analytic properties}

\subsection*{Relations with  $L$-functions and motives for $\mathrm{Sp}_n$
}
%
Following  \cite{Pa94} and \cite{Y}, these functions are conjectured to be
motivic for all $k>n$:
$$
   L(Sp(f), s, \chi) =L(M(Sp(f))(\chi), s), L(St(f), s) =L(M(St(f))(\chi), s), \mbox{ where}
$$
and the motives $M(Sp(f))$ and $M(St(f))$ are {\it pure} if $f$ is a genuine cusp form (not coming from a lifting of a smaller genus):
\begin{itemize}
{}\item
$M(Sp(f))$ is a motive over $\mathbb{Q}$ with coefficients in $\mathbb{Q}(\lambda_f(n))_{n\in\mathbb{N}}$
of  rank $2^n$, of weight $w=kn-n(n+1)/2$, and of  Hodge type
$\oplus_{p,q} H^{p,q}$, with
\begin{align}\label{Hodge} &
p=(k-i_1)+(k-i_2)+\cdots+(k-i_r),\\ \nonumber &
q=(k-j_1)+(k-j_2)+\cdots+(k-i_s),\mbox{ where }r+s=n,\\ \nonumber &
1\le i_1<i_2<\cdots<i_r\le n, 1\le j_1<j_2<\cdots<j_s\le n,\\ \nonumber &
\{ i_1,\cdots,i_r\}\cup \{ j_1,\cdots,i_s\}=\{ 1,2,\cdots,n\};
\end{align}
{}\item
$M(St(f))$ is a motive over $\mathbb{Q}$ with coefficients in $\mathbb{Q}(\lambda_f(n))_{n\in\mathbb{N}}$
of  rank $2n+1$, of weight $w=0$, and of  Hodge type
$H^{0,0}\oplus_{i=1}^n (H^{-k+i,k-i}\oplus H^{k-i, -k+i})$.
\end{itemize}
\subsection*{A functional equation 
}
Following  general Deligne's conjecture \cite{D}
on the motivic $L$-functions,
the $L$-function satisfy a functional equation  
determined by the Hodge structure of a motive:
$$
\Lambda(Sp(f), kn-n(n+1)/2+1-s)=  \varepsilon(f) \Lambda(Sp(f), s),  \text{~where}
$$
$
\Lambda(Sp(f), s)=\Gamma_{n,k}(s)L(Sp(f), s),  \varepsilon(f)=(-1)^{k2^{n-2}},
$ \\ \  \\ {}
 $\Gamma_{1,k}(s)=\Gamma_\mathbb{C}(s)=2(2\pi)^{-s}\Gamma(s)$, $\Gamma_{2,k}(s)=\Gamma_\mathbb{C}(s)\Gamma_\mathbb{C}(s-k+2)$, and\\
$\Gamma_{n,k}(s)=\prod_{p<q}\Gamma_\mathbb{C}(s-p)\Gamma^{a_+}_\mathbb{R}(s-(w/2))\Gamma_\mathbb{R}(s+1-(w/2))^{a_-}$ for some non-negative integers
$a_+$ and $a_-$, with
$a_++a_-=w/2$, and $\Gamma_\mathbb{R}(s)=\pi^{-s/2}\Gamma(s/2)$.

In particular, for $n=3 $ and $k\ge 5$, this conjectural functional equation has the form 
$\Lambda(Sp(f), s)=\Lambda(Sp(f), 3k-5-s)$, 
where
$$
\Lambda(Sp(f), s)=\Gamma_\CC(s)\Gamma_\CC(s-k+3)\Gamma_\CC(s-k+2)\Gamma_\CC(s-k+1)L(Sp(f), s).
$$
For $k\ge 5$ the critical values in the sense of Deligne \cite{D} are~:
$$
s=k, \cdots, 2k-5.
$$
\subsection*{A study of the analytic properties of $L(Sp(f), s)$
}
(compare with \cite{Vo}).

One could try to use a link between the eigenvalues  
$\lambda_f(T)$ 
and the 
Fourier coefficients $a_f({\mathcal T})$, where
$T\in {\mathcal D}(\Gamma, S)$ 
runs through the 
Hecke operators, and  ${\mathcal T}\in B_n $ 
runs over half-integral symmetric matrices. 
It was found by A.N.Andrianov (see \cite{An67}) that
\begin{align*}&
D({X})  = \sum_{\delta=0}^{\infty} \mathbf{T}(p^\delta ){X}^\delta 
=\frac{E({X})}{F({X})},
\end{align*}
where
\begin{align*}\label{P3}& \nonumber
E({X}) 
\\  \quad &
 =1-
p^2\left( \mathbf{T}_2(p^2)+
(p^2-p+1)(p^2+p+1)
[\mathbf{p}]_3\right){X}^2
+(p+1)p^4\mathbf{T}(p)[\mathbf{p}]_3{X}^3
\\ \nonumber \quad
      &
-p^7[\mathbf{p}]_3\left( \mathbf{T}_2(p^2)+
(p^2-p+1)(p^2+p+1)
[\mathbf{p}]_3\right) {X}^4 
   +p^{15}[\mathbf{p}]_3^3\,{X}^6 \, \in {\mathcal L}_\Z[{X}].
\end{align*}
\subsection*{Computing a formal Dirichlet series} 
Knowing $E(X)$, one computes the following  formal Dirichlet series
\begin{align*} &
D_E(s)  = \sum_{h=1}^{\infty} \mathbf{T}_E(h )h^{-s}=\prod_p D_{E,p}(p^{-s}),\mbox{ where}
\\ &
D_{E,p}(X)  = \sum_{\delta=0}^{\infty} \mathbf{T}_E(p^\delta )X^\delta 
=\frac{D_p(X)}{E(X)} 
=\frac{1}{F(X)}\in {\mathcal D}(\Gamma, S)[\![X]\!].
\end{align*}

\subsection*{Main identity}  
For all  ${\mathcal T}$ one obtains the following identity 
\begin{align*} &
a_f({\mathcal T})L(Sp(f), s)=\sum_{h=1}^{\infty} a_f(\mathcal{T}, E, h)h^{-s}, \mbox{ where}
\\ &
f|\mathbf{T}_E(h)
=\sum_{{\mathcal T}\in B_n} a_f(\mathcal{T}, E, h) q^{\mathcal T}.
\end{align*} 
Indeed, one has 
\begin{align*} &
f|\mathbf{T}_E(h)=\lambda_f(\mathbf{T}_E(h))f,\mbox{  and } L(Sp(f), s)=\sum_{h=1}^{\infty} \lambda_f(\mathbf{T}_E(h))h^{-s},
\mbox{   hence }
\\ &
 \sum_{h=1}^{\infty}f|\mathbf{T}_E(h)h^{-s}= L(Sp(f), s)\cdot f = \sum_{h=1}^{\infty}\sum_{{\mathcal T}\in B_n} 
a_f(\mathcal{T}, E, h) h^{-s} q^{\mathcal T},
\end{align*} 
and it suffices to compare the Fourier coefficients.
Such an identity is an unavoidable step in the problem of analytic continuation 
of $L(Sp(f), s)$  
and in the study of its arithmetical implications: one obtains a mean of computing special values out of 
Fourier coefficients.


For the standard  $L$-function $L(St(f), s)$, 
(see \cite{Pa94}, \cite{CourPa}, and 
\cite{Boe-Sch}).
\section{Critical values, periods and $p$-adic $L$-functions for $\Sp_3$
}
\subsection*{Critical values, periods and $p$-adic $L$-functions for $\Sp_3$}
A general conjecture by
Coates, Perrin-Riou
(see \cite{Co-PeRi}, \cite{Co}, \cite{Pa94}), 
 predicts that for $n=3$ and $k>5$, the motivic function $L(Sp(f), s)$, 
admits a $p$-adic  analogue. Let us fix an embedding 
$i_p: \overline {\mathbb Q}\to {\mathbb C}_p =
 \widehat{\overline {\mathbb Q}}_p$, and let $\alpha_0(p)$ be  an inverse root of $Q_{f,p}$ with the smallest $p$-valuation.
\begin{itemize}
\item
The {\it Panchishkin condition} (see \cite{Ha-Li-Sk}, \cite{Pa94}, \cite{PaAnnIF94}) for the existence of
  bounded $p$-adic $L$-functions in this case takes the form 
$\ord_p(\alpha_0(p))=0$  ($\alpha_0(p)$ is an inverse  root with the smallest $p$-valuation). 
Recall that  this condition  says:
\begin{quote} {\it {}
for a pure motive $M$ of  rank $d$, \\ Newton $p$-polygone at ($d/2$)=Hodge polygone at ($d/2$)}
\end{quote}
\item
Otherwise, one obtains $p$-adic $L$-functions of  logarithmic growth $o(\log^h(\cdot))$
with $h=[2\ord_p(\alpha_0(p))]+1$, $2\ord_p(\alpha_0(p))$ = the difference \\ 
Newton $p$-Polygone at ($d/2$)--Hodge Polygone at ($d/2$)
\end{itemize}
For the unitary groups, this 
 condition for the existence of $p$-adic  $L$-functions 
was discussed in   \cite{Ha-Li-Sk}. 

\section{Rankin's Lemma of higher genus}

Our next purpose  is to evaluate the generating series
$$
D^{(1,1)}_p(X)=\sum_{\delta=0}^\infty T(p^\delta)\otimes T(p^\delta)\,X^\delta\in\mathcal{L}_{2,\Z}\otimes\mathcal{L}_{2,\Z}[\![X]\!].
$$
in terms of Hecke algebra generators $\mathcal{L}_{2,\Z}\otimes\mathcal{L}_{2,\Z}$:
$$
\Ta\otimes 1,\quad \Tb\otimes 1,\quad \Te\otimes 1,
$$
$$
1\otimes\Ta,\quad 1\otimes\Tb,\quad 1\otimes\Te.
$$

\normalsize
\begin{theo}[\cite{PaVaRnk}]
\large
For genus $n=2$, we have the following explicit representation
\begin{align*}
&D^{(1,1)}_p(X)=\sum_{\delta =0}^\infty T(p^\delta)\otimes T(p^\delta)\,X^\delta=(1-p^6\Te\otimes\Te X^2)\frac{R(X)}{S(X)}\,,
\text{~where}\\\\
&R(X)=1+r_1X+\dotsc+r_{12}X^{12},\quad{\text{~with~~}r_{1}=r_{11}=0},\\
&S(X)=1+s_1X+\dotsc+s_{16}X^{16},\\
&\qquad R(X),S(X)\in\mathcal{L}_{2,\Z}\otimes\mathcal{L}_{2,\Z}[X],
\end{align*}
where $r_i$ and $s_i$ are given 
in the Appendix of \cite{PaVaRnk}
\end{theo}

\begin{figure}[h]
\caption{\label{NP} Newton polygons of $R(X)$ and $S(X)$ with respect to powers of $p$ and $X$,  of heights 34 and 48, resp.}
\begin{center}
\includegraphics[width=5cm,height=5cm]{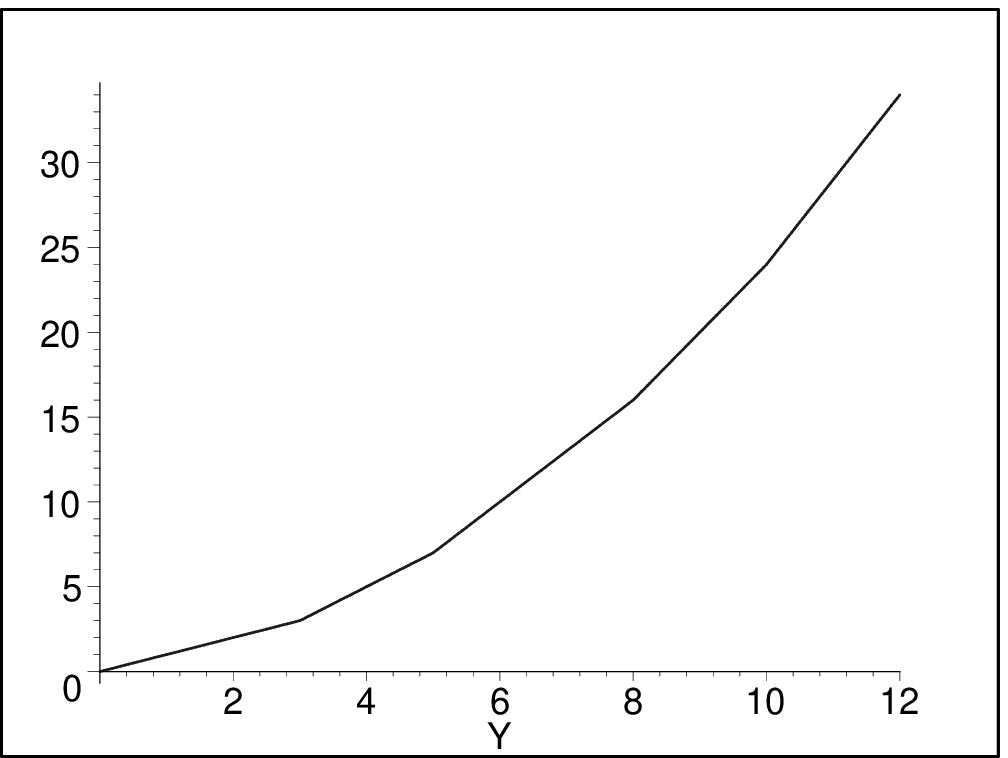}\hskip0.5cm
\includegraphics[width=5cm,height=5cm]{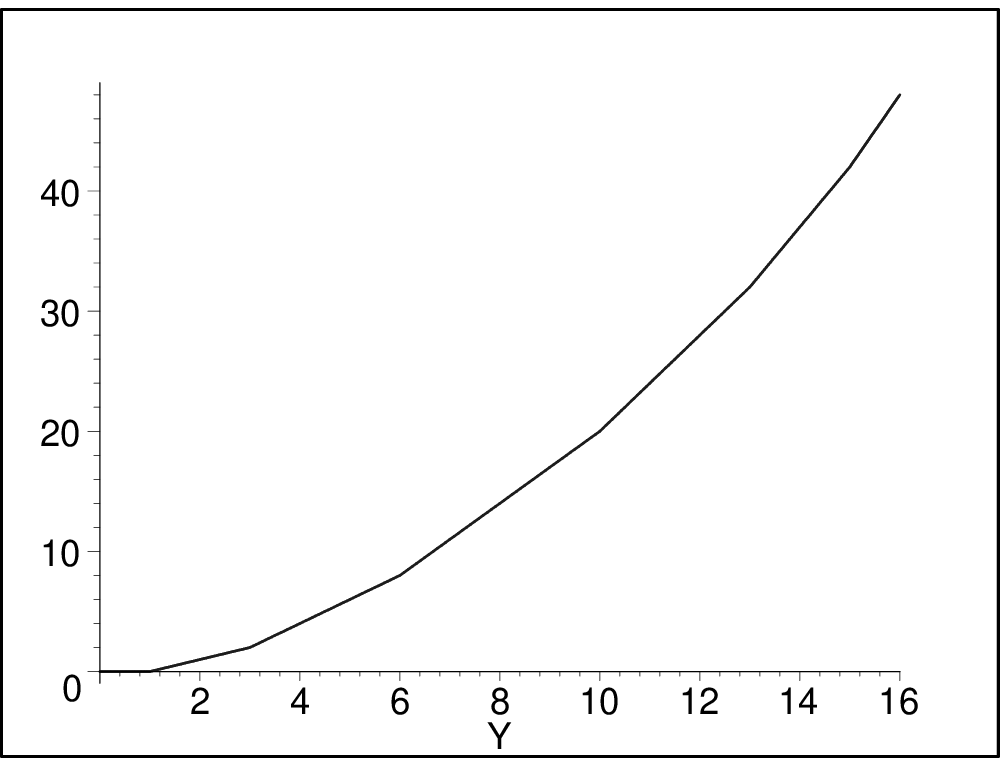}
\end{center}
\end{figure}

Here we give the Newton polygons of $R(X)$ and $S(X)$ with respect to powers of $p$ and $X$
(see Figure \ref{NP}). It follows from our comutation that all slopes are {\it integral}.
We hope that these polygons could  help to find some geometric objects
attached to
the polynomials  $R(X)$ and $S(X)$,
in the spirit of a recent work of C.Faber and G.Van Der Geer, \cite{FVdG}.

{
A half of  coefficients: $s_9,\dotsc,s_{16}$ can be found using this simple functional equation~:
$$
s_{16-i}=(p^6\Tc\otimes\Tc)^{8-i}s_i \ \ \ (i=0, \cdots, 8).
$$
}



{ 
{\it Proof: } 
Recall that
{ 
$$
\sum_{\delta=0}^\infty\,\Omega^{(n=2)}(T(p^\delta))X^\delta=\frac{1-\dfrac{x_0^2x_1x_2}{p}\,X^2}{(1-x_0X)\,(1-x_0x_1X)(1-x_0x_2X)\,(1-x_0x_1x_2X)\,}
$$
}
From this series it is possible to obtain a formula for ${\Omega }(T(p^\delta))$ considering four geometric progressions
{
\begin{align*}
&\sum_{\nu_1=0}^{\infty}(x_0X)^{\nu_1}= \frac{1}{1-x_0X},\quad\sum_{\nu_2=0}^{\infty}(x_0x_1X)^{\nu_2}=\frac{1}{1-x_0x_1X}\\ &\sum_{\nu_3=0}^{\infty}(x_0x_2X)^{\nu_3}=\frac{1}{1-x_0x_2X},\quad\sum_{\nu_4=0}^{\infty}(x_0x_1x_2X)^{\nu_4}=\frac{1}{1-x_0x_1x_2X}
\end{align*}
}
{
\begin{align*}
\Omega_x&(T(p^\delta))=p^{-1}x_0^\delta(px_1^{(3+\delta)}x_2-px_1x_2^{(3+\delta)}-px_1^{(2+\delta)}+px_2^{(2+\delta)}-px_1^{(3+\delta)}x_2^{(2+\delta)}\\
&+px_1^{(2+\delta)}x_2^{(3+\delta)}+px_1-px_2-x_1^{(2+\delta)}x_2^2+x_1^{(1+\delta)}x_2+x_1^{(2+\delta)}x_2^{(1+\delta)}-x_1^{(1+\delta)}x_2^{(2+\delta)}\\
&+x_1^2x_2^{(2+\delta)}-x_1x_2^{(1+\delta)}-x_1^2x_2+x_1x_2^2)/((1-x_1)(1-x_2)(1-x_1x_2)(x_1-x_2))\\
=&-p^{-1}x_0^\delta((1-x_1x_2)(px_1-x_2)x_1^{(\delta+1)}+(1-x_1x_2)(x_1-px_2)x_2^{(\delta+1)}\\
&-(1-px_1x_2)(x_1-x_2)(x_1x_2)^{(\delta+1)}-(p-x_1x_2)(x_1-x_2))/\\
&\quad((1-x_1)(1-x_2)(1-x_1x_2)(x_1-x_2)).
\end{align*}
}
}

{\subsubsection*{Tensor product of local Hecke algebras}
{ 
We use a second group of variables $y_0, y_1, y_2$ and $\Omega_y$ in order to obtain a tensor product of local Hecke algebras
\begin{align*}
&\Omega^{(n)}_x\otimes\Omega^{(n)}_y:\mathcal{L}_{n,\Q}\otimes\mathcal{L}_{n,\Q}\to\Q[x_0,x_1,\dotsc,x_n;y_0,y_1,\dotsc, y_n],\text{~et}\\
&\Omega^{(2)}_y(T(p^\delta)):=\Omega^{(2)}_x(T(p^\delta))|_{x=y}
\end{align*}
The product of two polynomials $\Omega^{(2)}_x(T(p^\delta))$ and $\Omega^{(2)}_y(T(p^\delta))$ is computed straightforward~:
\begin{align*}
\Omega&^{(2)}_x(T(p^\delta))\cdot\Omega^{(2)}_y(T(p^\delta))=p^{-2}x_0^\delta y_0^\delta
(px_1^{(3+\delta)}x_2-px_1^{(2+\delta)}-px_1^{(3+\delta)}x_2^{(2+\delta)}
\\&+px_1^{(2+\delta)}x_2^{(3+\delta)}-px_1x_2^{(3+\delta)}+px_2^{(2+\delta)}+px_1-px_2-x_1^{(2+\delta)}x_2^2+x_1^{(1+\delta)}x_2
\\&+x_1^{(2+\delta)}x_2^{(1+\delta)}-x_1^{(1+\delta)}x_2^{(2+\delta)}+x_1^2x_2^{(2+\delta)}-x_1x_2^{(1+\delta)}-x_1^2x_2+x_1x_2^2)
\\&\times(py_1^{(3+\delta)}y_2-py_1^{(2+\delta)}-py_1^{(3+\delta)}y_2^{(2+\delta)}+py_1^{(2+\delta)}y_2^{(3+\delta)}-py_1y_2^{(3+\delta)}\\&
+py_2^{(2+\delta)}+py_1-py_2-y_1^{(2+\delta)}y_2^2+y_1^{(1+\delta)}y_2+y_1^{(2+\delta)}y_2^{(1+\delta)}\\&
-y_1^{(1+\delta)}y_2^{(2+\delta)}+y_1^2y_2^{(2+\delta)}-y_1y_2^{(1+\delta)}-y_1^2y_2+y_1y_2^2)/
\\&\quad((1-x_1)(1-x_2)(1-x_1x_2)(x_1-x_2)(1-y_1)(1-y_2)(1-y_1y_2)(y_1-y_2))
\end{align*}
}
}

{\subsubsection*{}
Summation of obtained expression give the following result~:
{\footnotesize
\begin{align*}
&(\Omega^{(2)}\otimes\Omega^{(2)})(D^{(1,1)}_p(X)) =
\sum_{\delta=0}^\infty\,\Omega^{(2)}_x(T(p^\delta))\cdot\Omega^{(2)}_y(T(p^\delta))\,X^\delta =
\\ & \nonumber 
-{\displaystyle \frac 
 {(p\,\mathit{x}_1 - \mathit{x}_2)\,(1-p\,\mathit{y}_1\,\mathit{y}_2)\,
   \mathit{x}_1\,\mathit{y}_1\,\mathit{y}_2}
 {p^{2}\,(1-\mathit{x}_1)\,(1-\mathit{x}_2)\,(\mathit{x}_1-\mathit{x}_2)\,
   (1-\mathit{y}_1)\,(1-\mathit{y}_2)\,(1-\mathit{y}_1\,\mathit{y}_2)\,
   (1-\mathit{x}_0\,\mathit{x}_1\,\mathit{y}_0\,\mathit{y}_1\,\mathit{y}_2\,X)}}
\\ & \nonumber 
+{\displaystyle \frac 
{\mathit{x}_2\,\mathit{y}_1\,(\mathit{x}_1-p\,\mathit{x}_2)\,(p\,\mathit{y}_1-\mathit{y}_2)}
{p^{2}\,(1-\mathit{x}_1)\,(1-\mathit{x}_2)\,(\mathit{x}_1-\mathit{x}_2)\,
  (1-\mathit{y}_1)\,(1-\mathit{y}_2)\,(\mathit{y}_1-\mathit{y}_2)\,
  (1-\mathit{x}_0\,\mathit{x}_2\,\mathit{y}_0\,\mathit{y}_1\,X)}}
\\ & \nonumber 
+{\displaystyle \frac 
{\mathit{x}_2\,\mathit{y}_2\,(\mathit{x}_1-p\,\mathit{x}_2)(\mathit{y}_1-p\,\mathit{y}_2)}
{p^{2}\,(1-\mathit{x}_1)\,(1-\mathit{x}_2)\,(\mathit{x}_1-\mathit{x}_1)\,
  (1-\mathit{y}_1)\,(1-\mathit{y}_2)\,(\mathit{y1}-\mathit{y}_2)\,
  (1-\mathit{x}_0\,\mathit{y}_0\,\mathit{x}_2\,\mathit{y}_2\,X)}} 
\\ & \nonumber 
-{\displaystyle \frac 
{\mathit{x}_2\,\mathit{y}_1\,\mathit{y}_2\,(\mathit{x}_1-p\,\mathit{x}_2)\,
  (1-p\,\mathit{y}_1\,\mathit{y}_2)}
{p^{2}\,(1-\mathit{x}_1)\,(1-\mathit{x}_2)\,(\mathit{x}_1 -\mathit{x}_2)\,
  (1-\mathit{y}_1)\,(1-\mathit{y}_2)\,(1-\mathit{y}_1\,\mathit{y}_2)\,
(1-\mathit{x}_0\,\mathit{x}_2\,\mathit{y}_0\,\mathit{y}_1\,\mathit{y}_2\,X)}}
\\ & \nonumber 
-{\displaystyle \frac 
{\mathit{x}_1\,(p\,\mathit{x}_1-\mathit{x}_2)\,(p-\mathit{y}_1\,\mathit{y}_2)}
{p^{2}\,(1-\mathit{x}_1)\,(1-\mathit{x}_2)\,(\mathit{x}_1-\mathit{x}_2)\,
  (1-\mathit{y}_1)\,(1-\mathit{y}_2)\,(1-\mathit{y}_1\,\mathit{y}_2)\,
  (1-\mathit{x}_0\,\mathit{x}_1\,\mathit{y}_0\,X)}}
\\ & \nonumber 
-{\displaystyle \frac 
{\mathit{x}_1\,\mathit{x}_2\,\mathit{y}_1(1-p\,\mathit{x}_1\,\mathit{x}_2)\,
  (p\,\mathit{y}_1-\mathit{y}_2)}
{p^{2}\,(1-\mathit{x}_1)\,(1-\mathit{x}_2)\,(1-\mathit{x}_1\,\mathit{x}_2)\,
  (1-\mathit{y}_1)\,(1-\mathit{y}_2)\,(\mathit{y}_1-\mathit{y}_2)\,
  (1-\mathit{x}_0\,\mathit{x}_1\,\mathit{x}_2\,\mathit{y}_0\,\mathit{y}_1\,X)}}
\\ & \nonumber 
-{\displaystyle \frac 
{\mathit{x}_1\,\mathit{x}_2\,\mathit{y}_2\,(1-p\,\mathit{x}_1\,\mathit{x}_2)\,
  (\mathit{y}_1-p\,\mathit{y}_2)}
{p^{2}\,(1-\mathit{x}_1)\,(1-\mathit{x}_2)\,(1-\mathit{x}_1\,\mathit{x}_2)\,
  (1-\mathit{y}_1)\,(1-\mathit{y}_2)\,(\mathit{y}_1-\mathit{y}_2)\,
  (1-\mathit{x}_0\,\mathit{x}_1\,\mathit{x}_2\,\mathit{y}_0\,\mathit{y}_2\,X)}} \\ 
& \nonumber 
+{\displaystyle \frac 
{\mathit{y}_1\,\mathit{y}_2\,(p-\mathit{x}_1\,\mathit{x}_2)\,(1-p\,\mathit{y}_1\,\mathit{y}_2)}
{p^{2}\,(1-\mathit{x}_1)\,(1-\mathit{x}_2)\,(1-\mathit{x}_1\,\mathit{x}_2)\,
  (1-\mathit{y}_1)\,(1-\mathit{y}_2)\,(1-\mathit{y}_1\,\mathit{y}_2)\,
  (1-\mathit{x}_0\,\mathit{y}_0\,\mathit{y}_1\,\mathit{y}_2\,X)}}
\\ & \nonumber 
+{\displaystyle \frac 
{\mathit{x}_1\,\mathit{x}_2\,(1-p\,\mathit{x}_1\,\mathit{x}_2)\,
  (p-\mathit{y}_1\,\mathit{y}_2)}
{p^{2}\,(1-\mathit{x}_1)\,(1-\mathit{x}_2)\,(1-\mathit{x}_1\,\mathit{x}_2)\,
  (1-\mathit{y}_1)\,(1-\mathit{y}_2)\,(1-\mathit{y}_1\,\mathit{y}_2)\,
  (1-\mathit{x}_0\,\mathit{x}_1\,\mathit{x}_2\,\mathit{y}_0\,X)}}
\\ & \nonumber 
-{\displaystyle \frac 
{\mathit{x}_1\,\mathit{y}_1\,(p\,\mathit{x}_1-\mathit{x}_2)\,(p\,\mathit{y}_1-\mathit{y}_2)}
{p^{2}\,(1-\mathit{x}_1)\,(1-\mathit{x}_2)\,(\mathit{x}_1-\mathit{x}_2)\,
  (1-\mathit{y}_1)\,(1-\mathit{y}_2)\,(\mathit{y}_1-\mathit{y}_2)\,
  (1-\mathit{x}_0\,\mathit{x}_1\,\mathit{y}_0\,\mathit{y}_1\,X)}}
\\ & \nonumber 
+{\displaystyle \frac 
{\mathit{x}_1\,\mathit{y}_2\,(p\,\mathit{x}_1-\mathit{x}_2)\,(\mathit{y}_1-p\,\mathit{y}_2)}
{p^{2}\,(1-\mathit{x}_1)\,(1-\mathit{x}_2)\,(\mathit{x}_1-\mathit{x}_2)\,
  (1-\mathit{y}_1)\,(1-\mathit{y}_2)\,(\mathit{y}_1-\mathit{y}_2)\,
(1-\mathit{x}_0\,\mathit{x}_1\,\mathit{y}_0\,\mathit{y}_2\,X)}}
\\ & \nonumber 
-{\displaystyle \frac 
{\mathit{x}_2\,(\mathit{x}_1-p\,\mathit{x}_2)\,(p-\mathit{y}_1\,\mathit{y}_2)}
{p^{2}\,(1-\mathit{x}_1)\,(1-\mathit{x}_2)\,(\mathit{x}_1-\mathit{x}_2)\,
  (1-\mathit{y}_1)\,(1-\mathit{y}_2)\,(1-\mathit{y}_1\,\mathit{y}_2)\,
  (1-\mathit{x}_0\,\mathit{x}_2\,\mathit{y}_0\,X)}}
\\ & \nonumber 
+{\displaystyle \frac 
{\mathit{x}_1\,\mathit{x}_2\,\mathit{y}_1\,\mathit{y}_2\,(1-p\,\mathit{x}_1\,\mathit{x}_2)\,
  (1-p\,\mathit{y}_1\,\mathit{y}_2)}
{p^{2}\,(1-\mathit{x}_1)\,(1-\mathit{x}_2)\,(1-\mathit{x}_1\,\mathit{x}_2)\,
  (1-\mathit{y}_1)\,(1-\mathit{y}_2)\,(1-\mathit{y}_1\,\mathit{y}_2)\,
(1-\mathit{x}_0\,\mathit{x}_1\,\mathit{x}_2\,\mathit{y}_0\,\mathit{y}_1\,\mathit{y}_2\,X)}}
\\ & \nonumber 
+{\displaystyle \frac 
{(p-\mathit{x}_1\,\mathit{x}_2)\,(p-\mathit{y}_1\,\mathit{y}_2)}
{p^{2}\,(1-\mathit{x}_1)\,(1-\mathit{x}_2)\,(1-\mathit{x}_1\,\mathit{x}_2)\,
  (1-\mathit{y}_1)\,(1-\mathit{y}_2)\,(1-\mathit{y}_1\,\mathit{y}_2)\,
  (1-\mathit{x}_0\,\mathit{y}_0\,X)}}
\\ & \nonumber 
-{\displaystyle \frac 
{\mathit{y}_1\,(p-\mathit{x}_1\,\mathit{x}_2)\,(p\,\mathit{y}_1 - \mathit{y}_2)}
{p^{2}\,(1-\mathit{x}_1)\,(1-\mathit{x}_2)\,(1-\mathit{x}_1\,\mathit{x}_2)\,
  (1-\mathit{y}_1)\,(1-\mathit{y}_2)\,(\mathit{y}_1-\mathit{y}_2)\,
  (1-\mathit{x}_0\,\mathit{y}_0\,\mathit{y}_1\,X)}}
\\ & \nonumber 
-{\displaystyle \frac 
{\mathit{y}_2\,(p-\mathit{x}_1\,\mathit{x}_2)\,(\mathit{y}_1-p\,\mathit{y}_2)}
{p^{2}\,(1-\mathit{x}_1)\,(1-\mathit{x}_2)\,(1-\mathit{x}_1\,\mathit{x}_2)\,
  (1-\mathit{y}_1)\,(1-\mathit{y}_2)\,(\mathit{y}_1-\mathit{y}_2)\,
  (1-\mathit{x}_0\,\mathit{y}_0\,\mathit{y}_2\,X)}}\,.
\end{align*}
}

{\subsubsection*{Properties of the image $(\Omega^{(2)}\otimes\Omega^{(2)})(D^{(1,1)}_p(X))$}
We checked by an explicit computation that the polynomials, which do not depend of $X$ in the denominator of the image
 $(\Omega^{(2)}\otimes\Omega^{(2)})(D^{(1,1)}_p(X))$, disappear after simplification in the ring $\Q[x_0,x_1,x_2,y_0,y_1,y_2][\![X]\!]$, and the common denominator becomes~:
{ 
\begin{align*} &
(1 - \mathit{x}_0\,\mathit{y}_0 X)\,
(1 - \mathit{x}_0\,\mathit{y}_0\,\mathit{x}_1\,X)\,
(1 - \mathit{x}_0\,\mathit{y}_0\,\mathit{y}_1\,X)\,
(1 - \mathit{x}_0\,\mathit{y}_0\,\mathit{x}_2\,X)\,
(1 - \mathit{x}_0\,\mathit{y}_0\,\mathit{y}_2\,X)\,
 \\ &
(1 - \mathit{x}_0\,\mathit{y}_0\,\mathit{x}_1\,\mathit{y}_1\,X)\,
(1 - \mathit{x}_0\,\mathit{y}_0\,\mathit{x}_1\,\mathit{x}_2\,X)\,
(1 - \mathit{x}_0\,\mathit{y}_0\,\mathit{x}_1\,\mathit{y}_2\,X)\,
(1 - \mathit{x}_0\,\mathit{y}_0\,\mathit{y}_1\,\mathit{x}_2\,X)\,
 \\ &
(1 - \mathit{x}_0\,\mathit{y}_0\,\mathit{y}_1\,\mathit{y}_2\,X)\,
(1 - \mathit{x}_0\,\mathit{y}_0\,\mathit{x}_2\,\mathit{y}_2\,X)\,
(1 - \mathit{x}_0\,\mathit{y}_0\,\mathit{x}_1\,\mathit{y}_1\,\mathit{x}_2\,X)\,
(1 - \mathit{x}_0\,\mathit{y}_0\,\mathit{x}_1\,\mathit{y}_1\,\mathit{y}_2\,X)\,
 \\ &
(1 - \mathit{x}_0\,\mathit{y}_0\,\mathit{x}_1\,\mathit{x}_2\,\mathit{y}_2\,X)\,
(1 - \mathit{x}_0\,\mathit{y}_0\,\mathit{y}_1\,\mathit{x}_2\,\mathit{y}_2\,X)\,
(1 - \mathit{x}_0\,\mathit{y}_0\,\mathit{x}_1\,\mathit{y}_1\,\mathit{x}_2\,\mathit{y}_2\, X).
\end{align*}
}
Moreover, we found that the numerator consists of the factor $(1-x_0^2y_0^2x_1y_1x_2y_2X^2)$ and a polynomial in $X$ of degree $12$ with coefficients in $\Q[x_0,x_1,x_2,y_0,y_1,y_2]$
 (the constant term is equal to $1$ and the main term is $p^{-2}x_0^{12}y_0^{12}x_1^6x_2^6y_1^6y_2^6X^{12}$). 
We also found that the factor of degree $12$ does not contain terms of degree $1$ and $11$ in $X$.

It follows that {
$(\Omega^{(2)}\otimes\Omega^{(2)})(D^{(1,1)}_p(X))={
\dfrac{(1-x_0^2y_0^2x_1y_1x_2y_2X^2)R_{x,y}(X)}{S_{x,y}(X)}}$}, where $R_{x,y}(X)=1+r_{2,x,y}X^2+\dotsb+r_{10,x,y}X^{10}
+r_{12,x,y}X^{12}$ and $S_{x,y}(X)=1+s_{1,x,y}X+\dotsb+s_{16,x,y}X^{16}\in\Q[x_0,x_1,x_2,y_0,y_1,y_2,X]$.
}

{\subsubsection*{Expression through the Hecke operators}
Knowing the coefficients of $R_{x,y}(X)$, $S_{x,y}(X)\in\Q[x_0,x_1,x_2,y_0,y_1,y_2,X]$
 (the images under $\Omega_x\otimes\Omega_y$)
 one can reconstruct $R(X)=1+r_2X^2+\dotsb+r_{10}X^{10}+r_{12}X^{12}$ and $S(X)=1+s_1X
+\dotsb+s_{16}X^{16}\in\mathcal{L}_{2,\Z}\otimes\mathcal{L}_{2,\Z}[X]$.

\bigbreak
Using the images of the generators
$$
\Omega_x(\Ta^{\lambda_0}\Tb^{\lambda_1}\Te^{\lambda_2})\Omega_y(\Ta^{\mu_0}\Tb^{\mu_1}\Te^{\mu_2}),
$$
we build and solve a linear system of undetermined coefficients $K_{\lambda_0,\lambda_1,\lambda_2,\mu_0,\mu_1,\mu_2}$
 expressing all monomials
 $\Omega_x(\Ta^{\lambda_0}\Tb^{\lambda_1}\Te^{\lambda_2})\Omega_y(\Ta^{\mu_0}\Tb^{\mu_1}\Te^{\mu_2})$,
 up to degree $12$ for $R_{x,y}(X)$ ($\leqslant 16$ for $S_{x,y}(X)$) in $x_0$ and in $y_0$.\qed
}

{\subsection*{Comparison with the case of genus $n=1$ 
}
The factor $(1-\mathit{x}_0^{2}\mathit{y}_0^{2}\,\mathit{x}_1\mathit{y}_1\mathit{x}_2\mathit{y}_2X^{2})$ of degree $2$ in $X$
 is very similar to the one in case of $g=1$~:
{
\begin{align*}
 & \sum _{\delta =0}^{\infty }\,
{\Omega ^{(1)}_{x}}(T(p^{\delta }))\cdot 
{\Omega ^{(1)}_{y}}(T(p^{\delta }))\,X^{\delta } = 
\sum _{\delta =0}^{\infty }
{\displaystyle \frac 
{\mathit{x}_0^{\delta }\, (1-\mathit{x}_1^{(1+\delta )})}{1-\mathit{x}_1}} 
\cdot
{\displaystyle \frac 
{\mathit{y}_0^{\delta }\, (1-\mathit{y}_1^{(1+\delta )})}{1-\mathit{y}_1}} 
X^{\delta }
\\ &  =
{\displaystyle \frac {1}
{(1-\mathit{x}_1)\,(1-\mathit{y}_1)\,(1-\mathit{x}_0\,\mathit{y}_0\,X)}}
 - {\displaystyle \frac {\mathit{y}_1}
{(1-\mathit{x}_1)\,(1-\mathit{y}_1)\,(1-\mathit{x}_0\,\mathit{y}_0\,\mathit{y}_1\,X)}} 
 \\ &
- {\displaystyle \frac {\mathit{x}_1}
{(1-\mathit{x}_1)\,(1-\mathit{y}_1)\,(1-\mathit{x}_0\,\mathit{y}_0\,\mathit{x}_1\,X)}}
+ {\displaystyle \frac {\mathit{x}_1\,\mathit{y}_1}
{(1-\mathit{x}_1)\,(1-\mathit{y}_1)\,(1-\mathit{x}_0\,\mathit{y}_0\,\mathit{x}_1\,\mathit{y}_1\,X)}}
\\ & \\ & =
{
\displaystyle \frac {1-\mathit{x}_0^{2}\,\mathit{y}_0^{2}\,\mathit{x}_1\,\mathit{y}_1\,X^{2}}
{(1-\mathit{x}_0\,\mathit{y}_0\,X) 
 (1-\mathit{x}_0\,\mathit{y}_0\,\mathit{x}_1\,X)\,
 (1-\mathit{x}_0\,\mathit{y}_0\,\mathit{y}_1\,X)\,
 {(1-\mathit{x}_0\,\mathit{y}_0\,\mathit{x}_1\,\mathit{y}_1\,X)}}}.
\end{align*}
\begin{align*}
\sum_{\delta=0}^\infty&\,T(p^{\delta })\otimes T(p^{\delta })\,X^\delta
=(1-p^2\Te\otimes\Te X^{2})\times\\
&\times(1-\Ta\otimes\Ta X+(p\Ta^2\otimes\Te+p\Te\otimes\Ta^2-2p^2\Te\otimes\Te)X^2\\
&-p^2\Ta\Te\otimes\Ta\Te X^3+p^4\Te^2\otimes\Te^2 X^4)^{-1}\,.
\end{align*}
}
}

\section{A lifting from $\,GSp_{\,2m} \times GSp_{\,2m}$ to $GSp_{\,4m}$ of genus $g=4m$} 
\subsection*{Motive of the  Rankin product of genus $n=2$}
Let $f$ and $g$ be two  Siegel cusp eigenforms of weights $k$ and $l$, $k>l$, and
let $M(Sp(f))$ and $M(Sp(g))$ be the (conjectural) spinor motives of  $f$ and $g$.
Then  $M(Sp(f))$ is a motive over $\mathbb{Q}$ with coefficients in $\mathbb{Q}(\lambda_f(n))_{n\in\mathbb{N}}$
of rank  $4$, of weight $w=2k-3$, and of  Hodge type
$H^{0,2k-3}\oplus H^{k-2,k-1}\oplus H^{k-1,k-2}\oplus H^{2k-3,0}$,
and
$M(Sp(g))$ is a motive over $\mathbb{Q}$ with coefficients in $\mathbb{Q}(\lambda_g(n))_{n\in\mathbb{N}}$
of rank $4$, of weight $w=2l-3$, and of  Hodge type
$H^{0,2l-3}\oplus H^{l-2,l-1}\oplus H^{l-1,l-2}\oplus H^{2l-3,0}$.

For another Siegel modular form, an eigenfunction of Hecke operators
  $g \in\mathcal{M}_l^n$ consider the corresponding homomorphism
$\lambda_g : {\mathcal L}_{n,p} \longrightarrow \mathbb{C}$ 
 given by its Satake parameters  $ (\beta_0, \beta_1, \cdots, \beta_n)$ of $g$, and let
$\lambda_f\otimes\lambda_g : 
{\mathcal L}_{n,p}\otimes{\mathcal L}_{n,p} \longrightarrow \mathbb{C}$


\subsection*{The tensor product $M(Sp(f))\otimes M(Sp(g))$}
is a motive over $\mathbb{Q}$ with coefficients in  $\mathbb{Q}(\lambda_f(n),\lambda_g(n))_{n\in\mathbb{N}}$
of rank $16$, of weight $w=2k+2l-6$, and of  Hodge type

\begin{align*}&
H^{0,2k+2l-6}\oplus H^{l-2,2k+l-4}\oplus H^{l-1,2k+l-5}\oplus H^{2l-3,2k-3}
\\ &
H^{k-2,k+2l-4}\oplus H^{k+l-4,k+l-2}\oplus H_+^{k+l-3,k+l-3}\oplus H^{k+2l-5,k-1}
\\ &
H^{k-1,k+2l-5}\oplus H_-^{k+l-3,k+l-3}\oplus H^{k+l-2,k+l-4}\oplus H^{k+2l-4,k-2}
\\ &
H^{2k-3,2l-3}\oplus H^{2k+l-5,l-1}\oplus H^{2k+l-4,l-2}\oplus H^{2k+2l-6,0}.
\end{align*}
\subsection*{Motivic  $L$-functions: analytic properties
}
Following  Deligne's conjecture \cite{D}
on motivic $L$-functions, applied
for a Siegel cusp eigenform $F$ for the Siegel modular group $\mathrm{Sp}_4(\mathbb{Z})$ of genus $n=4$ and  of weight  $k>5$,
one has
$\Lambda(Sp(F), s)=\Lambda(Sp(F), 4k-9-s)$,  
where
\begin{align*}&
\Lambda(Sp(F), s)=\Gamma_\mathbb{C}(s)\Gamma_\mathbb{C}(s-k+4)\Gamma_\mathbb{C}(s-k+3)\Gamma_\mathbb{C}(s-k+2)\Gamma_\mathbb{C}(s-k+1)
\\ & \times
\Gamma_\mathbb{C}(s-2k+7)\Gamma_\mathbb{C}(s-2k+6)\Gamma_\mathbb{C}(s-2k+5)
L(Sp(F), s),
\end{align*}
(compare this functional equation with that given in  \cite{An74}, p.115).

\

On the other hand,
for $m=2$ and for two cusp eigenforms  $f$ and $g$ for $\mathrm{Sp}_2(\mathbb{Z})$  of weights $k,l$, $k>l+1$,
$\Lambda(Sp(f)\otimes Sp(g), s)=\varepsilon(f,g)\Lambda(Sp(f)\otimes Sp(g), 2k+2l-5-s)$, $|\varepsilon(f,g)|=1$, where
\begin{align*}&
\Lambda(Sp(f)\otimes Sp(g), s)=\Gamma_\mathbb{C}(s)\Gamma_\mathbb{C}(s-l+2)\Gamma_\mathbb{C}(s-l+1)\Gamma_\mathbb{C}(s-k+2)
\\ & \times
\Gamma_\mathbb{C}(s-k+1)\Gamma_\mathbb{C}(s-2l+3)\Gamma_\mathbb{C}(s-k-l+2)\Gamma_\mathbb{C}(s-k-l+3)
\\ & \times
L(Sp(f)\otimes Sp(g), s).
\end{align*}
We used here the Gauss duplication formula $\Gamma_\mathbb{C}(s)=\Gamma_\mathbb{R}(s)\Gamma_\mathbb{R}(s+1)$.
Notice that  $a_+=a_-=1$ in this case,
and the conjectural motive $M(Sp(f))\otimes M(Sp(g))$ does not admit critical values.

\subsection*{A holomorphic lifting from $GSp_{2m} \times GSp_{2m}$ to $GSp_{4m}$: 
a conjecture}
\begin{conj}[on a lifting from $\,GSp_{\,2m} \times GSp_{\,2m}$ to $GSp_{\,4m}$]\label{lift2m}

Let $f$ and $g$ be two Siegel modular forms of {\it genus} $2m$ and of weights $k> 2m$ and $l=k-2m$.
Then there exists a Siegel modular form $F$ of genus $4m$
 and of weight $k$ with the Satake parameters
$\gamma_0=\alpha_0\beta_0, \gamma_1=\alpha_1, 
\gamma_2=\alpha_2, \cdots, \gamma_{2m}=\alpha_{2m},
\gamma_{2m+1}=\beta_1, \cdots, \gamma_{4m}=\beta_{2m}$
 for suitable choices
$\alpha_0, \alpha_1, \cdots, \alpha_{2m}$ 
and  $\beta_0, \beta_1, \cdots, \beta_{2m}$
of Satake's parameters of $f$ and $g$.


One readily checks that the Hodge types of $M(Sp(f))\otimes M(Sp(g))$  and $M(Sp(F))$ are the same (of rank $2^{4m}$)
(it follows from the above description (\ref{Hodge}), and from Künneth's-type formulas).

\end{conj}
An evidence for this version of the conjecture comes from Ikeda-Miyawaki constructions (\cite{Ike01}, \cite{Ike06}, \cite{Mur02}):
let $k$ be an even positive integer,
${h}\in S_{2k}(\Gamma_1)$ a normalized Hecke eigenform of weight $2k$,
$F_{2n}\in S_{k+n}(\Gamma_{2n})$ the Ikeda lift of ${h}$ of genus $2n$
(we assume $k\equiv n\bmod 2$, $n\in \mathbb{N}$).

Next let
$f\in S_{k+n+r}(\Gamma_{r})$  be an arbitrary   Siegel cusp eigenform of genus $r$ and weight $k+n+r$,  with $n,r\ge 1$.
If we take 
$n=m, r=2m$, $k:=k+m$, $k+n+r:=k+3m$,
then an example of the validity of this version of the conjecture is
 given by
\begin{align*} &
(f, g)=(f, F_{2m}(h))\mapsto {\cal F}_{{h},f}\in S_{k+3m}(\Gamma_{4m}), \\ &
(f, g)=(f, F_{2m})\in   S_{k+3m}(\Gamma_{2m})\times S_{k+m}(\Gamma_{2m}).
\end{align*}
\subsection*{}
Another evidence comes from Siegel-Eiseinstein series
$$
f=E^{2m}_k  \mbox{ and  }g=E^{2m}_{k-2m}
$$
 of  even genus $2m$ and weights $k$ and $k-2m$:
we have then
\begin{align*}&
\alpha_0=1, \alpha_1=p^{k-2m}, \cdots, \alpha_{2m}=p^{k-1},
\\ &
\beta_0=1, \beta_1=p^{k-4m}, \cdots, \beta_{2m}=p^{k-2m-1},
\end{align*}
then we have that
$$
\gamma_0=1, \gamma_1=p^{k-4m}, \cdots, \gamma_{2m}=p^{k-1},
$$
are the Satake parameters of the Siegel-Eisenstein series $F=E^{4m}_{k}$.
\subsection*{}
\begin{rema}
If we compare the $L$-function of the conjecture
(given by  the Satake parameters
$\gamma_0=\alpha_0\beta_0, \gamma_1=\alpha_1, 
\gamma_2=\alpha_2, \cdots, \gamma_{2m}=\alpha_{2m},
\gamma_{2m+1}=\beta_1, \cdots, \gamma_{4m}=\beta_{2m}$
 for suitable choices
$\alpha_0, \alpha_1, \cdots, \alpha_{2m}$ 
and  $\beta_0, \beta_1, \cdots, \beta_{2m}$
of Satake's parameters of $f$ and $g$), 
we see that it 
corresponds to the tensor product of the spinor $L$-functions, 
and this function is {\it not of the same type}
as that of the Yoshida's lifting \cite{Y81},  
which is a certain product of Hecke's $L$-functions.
\end{rema}
We would like 
to mention in this context Langlands's  functoriality: 
The denominators of our $L$-series
belong to local Langlands $L$-factors (attached to representations of $L$-groups). 
If we consider the homomorphisms 
$$
{}^LGSp_{2m}=GSpin(4m+1)\to 
GL_{2^{2m}},  \ \ 
{}^LGSp_{4m}=GSpin(8m+1)\to 
GL_{2^{4m}} ,
$$
we see that our conjecture is compatible 
with the homomorphism of 
$L$-groups
$$
GL_{2^{2m}} \times GL_{2^{2m}} \to 
GL_{2^{4m}} , \ (g_1, g_2)\mapsto g_1\otimes g_2, 
GL_n({\mathbb C})= {}^LGL_n.
$$
 However, it is unclear to us if Langlands's functoriality
 predicts a holomorphic Siegel modular form as a lift. 
\section{Constructions of  $p$-adic families of Siegel modular forms}
\subsection*{Constructing $p$-adic $L$-functions and modular symbols}
Together with complex parameter $s$ it is possible to use certain $p$-adic parameters in order to construct $L$-functions and modular symbols.  
We study such parameters as the twist with Dirichlet character on one hand, and the weight parameter in the theory
 of families of
 modular forms on the other.  
An operation of the twist by a Dirichlet character is a fundamental operation concerning the formal power series,
 and the $p$-adic variation of such character gives an analytic family of modular forms.  
One computes the modular symbols related to an automorphic representation $\pi$ of an algebraic group $G$
 over a number field,
 using the twists of automorphic $L$-functions $L(s,\pi\otimes\chi,r)$ with Dirichlet character $\chi$.  
We represent their special values as integrals giving both complex-analytic and $p$-adic-analytic continuation.  
We construct analytic families of such $L$-functions in cases $G=\mathrm{GL}_2\times\mathrm{GL}_2\times\mathrm{GL}_2$,
 $G=\mathrm{GL}_2\times\mathrm{GSp}_{2m}$, $G=\mathrm{GSp}_{2m}\times\mathrm{GSp}_{2m}$
 using the doubling method and its $p$-adic versions, which expected to work also for overconvergent families
 of automorphic forms,
 and already developed in simpler situations $G=\mathrm{GL}_2$ by R.~Coleman, G.~Stevens and others.  
\subsection*{A $p$-adic approach}
Consider Tate's field $\mathbb{C}_p={\widehat{\overline{\mathbb{Q}}}}_p$ for prime number $p$.  
Let us fix the embedding ${\overline{\mathbb{Q}}}\buildrel i_p\over\hookrightarrow\mathbb{C}_p$ and
 consider the algebraic numbers as numbers $p$-adic over $i_p$.  
For a $p$-adic family $k{}\mapsto f_k=\sum_{n=1}^\infty a_{n}(k)q^n \in {\overline{\mathbb{Q}}}[\![q]\!] \subset \mathbb{C}_p[\![q]\!]$, 
the  
Fourier coefficients $a_n(k)$ of $f_k$ and one of  Satake $p$-parameters  $\alpha(k):=\alpha_p^{(1)}(k)$ are given by the 
certain analytical $p$-adic functions $k{}\mapsto a_n(k)$ for $(n,p)=1$.  
The archetypal example of $p$-adic family is given by the Eisenstein series.
$$
a_n(k) = \sum_{d|n{},{}{}(d,p)=1} d^{k-1}, {}{}f_k=E_k, \alpha_p^{(1)}(k)=1, \alpha_p^{(2)}(k)=p^{k-1}.
$$
\noindent
The existence of the $p$-adic family of cusp forms of positive slope $\sigma >0$ was shown by Coleman. 
We define {\em the slope} $\sigma= {\mathrm{ord}}_p(\alpha_p^{(1)}(k))$ (and ask it to be constant in a $p$-adic neighborhood of weight $k$)
An example for $p=7$, $f=\Delta$, $k=12$, $a_7=\tau(7)=-7\cdot2392$, $\sigma=1$ is given by R.~Coleman in \cite{CoPB}.
\subsection*{Motivation for considering of $p$-adic families}
 comes from the conjecture of Birch and Swinnerton-Dyer, see \cite{Colm03}.  
For a cusp eigenform $f=f_2$, corresponding to an elliptic curve $E$ by Wiles \cite{Wi95}, we consider a family containing $f$.  
One can try to approach $k=2, s=1$ from the direction, taking $k\to 2$, instead of $s\to 1$, this leads to a formula linking the derivative over $s$
 at $s=1$ of the $p$-adic $L$-function with the derivative over $k$ at $k=2$ of the $p$-adic analytic function $\alpha_p(k)$,
 see in \cite{CST98}:
 $
{L'_{p,f}(1)=\mathcal{L}_p(f)L_{p,f}(1)}$ with $\mathcal{L}_p(f) = -2 \dfrac{d\alpha_p(k)}{dk} \big|_{k=2}$.  
The validity of this formula needs the existence of our two variable $L$-function,   constructed in \cite{PaTV}.

In order to construct $p$-adic $L$-function of two variables $(k,s)$, the theory of $p$-adic integration is used.  
The $H$-admissible measures with integer $H$ related to $\sigma$, which appear in this construction, 
are obtained from $H$-admissible measures with values in various rings of modular forms, in particular nearly holomorphic modular forms.  
Let us describe a typical example related to differential operators.
\subsection*{Nearly holomorphic modular forms and the method of canonical projection }
Let $\mathcal{A}$ be a commutative field.  There are several $p$-adic approaches for studying special values of $L$-functions using
 the method of canonical projection (see \cite{PaTV}).  
In this method the special values and modular symbols are considered as $\mathcal{A}$-linear forms over spaces of modular forms
 with coefficients in $\mathcal{A}$.  Nearly holomorphic (\cite{ShiAr}) modular forms are certain formal series
$$
g=\sum_{n=0}^{\infty}a(n;R)q^n \in \mathcal{A}[\![q]\!][R]
$$
with the property for $\mathcal{A}=\mathbb{C}$, $z=x+iy\in\mathbb{H}$, $R=(4\pi y)^{-1}$,
 the series converge to a ${\mathcal{C}}^{\infty}$-modular form over $\mathbb{H}$ of given weight $k$ and Dirichlet character $\psi$.  
The coefficients $a(n;R)$ are polynomials in $\mathcal{A}[R]$ of bounded degree.  
\subsection*{Triple products families} 
give a recent example of families on the algebraic group of higher rank.  
This aspect was studied by S.~Boecherer and A.~Panchishkin \cite{Boe-Pa6}.  
The triple product with Dirichlet character $\chi$ is defined as a complex $L$-function (Euler product of degree 8)
 $$
L(f_1\otimes f_2\otimes f_3, s, \chi) = \prod_{p\nmid N} L((f_1\otimes f_2\otimes f_3)_p, \chi(p) p^{-s}),
$$
 where
\begin{equation*}
\begin{split}
L((f_1\otimes f_2\otimes f_3)_p, X)^{-1}&=\det
\left(1_8-X\begin{pmatrix}\alpha^{(1)}_{p,1}&0\\0&\alpha^{(2)}_{p,1}\end{pmatrix}
    \otimes\begin{pmatrix}\alpha^{(1)}_{p,2}&0\\0&\alpha^{(2)}_{p,2}\end{pmatrix}
    \otimes\begin{pmatrix}\alpha^{(1)}_{p,3}&0\\0&\alpha^{(2)}_{p,3}\end{pmatrix}
\right).
\end{split}
\end{equation*}

We use the corresponding normalized $L$-function (see \cite{D}, \cite{Co}, \cite{Co-PeRi})~: 
$$
\Lambda (f_1\otimes f_2\otimes f_3, s, \chi) = 
\Gamma_{\mathbb{C}}(s) \Gamma_{\mathbb{C}}(s-k_3+1)
\Gamma_{\mathbb{C}}(s-k_2+1)\Gamma_{\mathbb{C}}(s-k_1+1) 
L(f_1\otimes f_2\otimes f_3, s, \chi),
$$
 where $\Gamma_{\mathbb{C}}(s) = 2(2\pi )^{-s}\Gamma(s)$.  
The Gamma-factor determines the {\em critical values} $s=k_1,\dotsc,k_2+k_3-2$ of $\Lambda(s)$, 
which we explicitly evaluate (like in the classical formula $\zeta(2)=\frac{\pi^2}{6}$). 
{\em A  functional equation} of $\Lambda(s)$ has the type 
$$
s\mapsto k_1+k_2+k_3-2-s.
$$

Let us consider the product of three eigenvalues: $\lambda= \lambda(k_1{}, k_2{}, k_3{})
= \alpha_{p, 1}^{(1)}(k{}_1)\alpha_{p, 2}^{(1)}(k{}_2)\alpha_{p, 3}^{(1)}(k_3{})$
with the slope 
$\sigma= v_p(\lambda(k_1{}, k_2{}, k_3{}))=
\sigma(k_1{}, k_2{}, k_3{})=
\sigma_1+\sigma_2+\sigma_3$
  {\it constant and positive} for all triplets $(k_1{}, k_2{}, k_3{})$ 
 in an appropriate $p$-adic neighbourhood
 of the fixed triplet of weights $(k_1, k_2, k_3)$.
\subsection*{The statement of the problem} 
for triple products is the following:
 {\em given three $p$-adic analytic families $\textbf{f}_j$ of slope $\sigma_j\ge0$,
 to construct a four-variable $p$-adic $L$-function attached to Garrett's triple product of these families}.  
We show that this function interpolates the special values 
$(s,k_1,k_2,k_2)\longmapsto\Lambda(f_{1,k_1}\otimes f_{2,k_2}\otimes f_{3,k_3},s,\chi)$
 at critical points $s=k_1,\dotsc,k_2+k_3-2$ for balanced weights $k_1 \le k_2+k_3-2$; 
we prove that these values are algebraic numbers after dividing by certain ``periods''.  
However the construction uses directly modular forms, and not the $L$-values in question, 
and a comparison of special values of two functions is done {\em after the construction}.  

\subsection*{Main result for triple products}
\begin{itemize}
\item[1)] The function 
$\mathcal{L}_{\textbf{f}}:(s,k_1,k_2,k_3)\mapsto\frac{\langle{\textbf{f}}^0,
\mathcal{E}(-r,\chi)\rangle}{\langle{\textbf{f}}^0,{\textbf{f}}_0\rangle}$
depends $p$-adic analytically on four variables \\ $(\chi\cdot y_p^r,k_1,k_2,k_3)\in X\times{\mathcal{B}}_1\times{\mathcal{B}}_2\times{\mathcal{B}}_3$;
\item[2)] Comparison of complex and $p$-adic values: for all $({k_1, k_2, k_3})$ in an affinoid neighborhood \\ ${\mathcal{B}}={\mathcal{B}}_1\times{\mathcal{B}}_2\times{\mathcal{B}}_3\subset X^3$, satisfying $k{}_1 \le k_2 + k_3 -2$: the values at $s=k_2 + k_3 -2-r$ \\ coincide with the normalized critical special values
 $$
L^*({f_{1,k_1}}\otimes{f_{2,k_2}}\otimes{f_{3,k_3}},k_2+k_3-2-r,\chi) \ \ (r=0,\dotsc,k_2+k_3-k_1-2)\,,
$$
for Dirichlet characters $\chi\bmod Np^v, v\ge 1$.
\item[3)] Dependence on $x\in X$: let $H= [2{\rm ord}_p(\lambda)]+1$.  
For any fixed $({k_1, k_2, k_3})\in\mathcal{B}$ and $x=\chi\cdot y_p^r$ then the following linear form 
(representing modular symbols for triple modular forms)
$$
x\longmapsto\frac{\left\langle\textbf{f}^0,\mathcal{E}(-r,\chi)\right\rangle}{\left\langle\textbf{f}^0,\textbf{f}_0\right\rangle,}
$$
extends to a $p$-adic analytic function of type $o(\log^H(\cdot))$ of the variable $x\in X$.
\end{itemize}

\subsection*{A general program}

We plan to extend this construction to other sutuations as follows~:

\begin{itemize}

\item[1)] Construction of modular distributions $\Phi_j$ 
with values in an infinite dimensional modular tower $\mathcal{M}(\psi)$.

\item[2)] Application of a canonical projector  of type $\pi_\alpha$ onto a finite dimensional subspace 
${\mathcal{M}}^{\alpha}(\psi)$ of ${\mathcal{M}}^{\alpha}(\psi)$. 

\item[3)] General admissibility criterium. 
The family of distributions $\pi_\alpha(\Phi_j)$ 
with value in $\mathcal{M}^\alpha(\psi)$ 
give a $h$-admissible measure $\tilde\Phi$
 with value in a module of finite rank.

\item[4)] Application of a linear form $\ell$ of type of a modular symbol
produces distributions $\mu_j=\ell (\pi_\alpha(\Phi_j))$, 
and  an admissible measure
from congruences between modular forms $\pi_\alpha(\Phi_j)$.

\item[5)] One shows that certain integrals $\mu_j(\chi)$ of the distributions $\mu_j$ coincide with certain 
$L$-values; however, these integrals are not necessary for the construction of measures 
(already done at stage 4).

\item[6)] One shows a result of uniqueness for the constructed $h$-admissible measures~:
they are determined by many of their integrals over Dirichlet characters (not all). 

\item[7)] In most cases we can prove a functional equation for the constructed measure $\mu$
 (using the uniqueness in 6),
 and using a functional equation for the $L$-values (over complex numbers, computed at stage 5).

\end{itemize}

\noindent
This strategy is already applicable in various cases
\section{Ikeda-Miyawaki constructions and their $p$-adic versions}
\subsection*{Ikeda's lift}
Ikeda generalized in 1999  (see \cite{Ike01}) the Saito-Kurokawa
lift of modular forms from one variable to Siegel modular forms of degree 2: 
 under the condition that ${n} \equiv k (\bmod2)$ there exists a lifting from an eigenform
$f\in   S_{2k}(\Gamma_1)$ to an eigenform $F\in S_{{n}+k}(\Gamma_{2{n}})$ such that the standard zeta function
$L(St(F), s)$ of $F$ (of degree $2{n}$) is given in terms of that of $f$ by
$$
\zeta(s)
\prod^{2{n}}
_{j=1}
L(f; s + k + {n}- j).
$$
(it was conjectured by Duke and Imamo\v glu in  \cite{DI98}). 
Notice that the Satake parameters of $F$ can be chosen in the form 
$
\beta_0, \beta_1, \dots, \beta_{2{n}}$, where 
$$
\beta_0=p^{{n}k-{n}({n}+1)/2},
\beta_i=\tilde\alpha p^{i-1/2} (i=1, \dots, {n}), \beta_{{n}+i}=\tilde\alpha^{-1} p^{i-1/2},
$$ 
 and $(1-\tilde\alpha p^{k-1/2}X)(1-\tilde\alpha^{-1}p^{k-1/2}X)=1-a(p)X+p^{2k-1}X^2$, see \cite{Mur02}, Lemma 4.1, p.65
  (so that  $\alpha=\tilde\alpha p^{k-1/2}$, $\tilde\alpha=\alpha p^{1/2-k}$ in our previous  notation).

\subsection*{Ikeda-Miyawaki conjecture}
(it was conjectured by Miyawaki in \cite{Mi92} and proved by Ikeda in \cite{Ike06}).

Let $k$ be an even positive integer, 
$f\in S_{2k}(\Gamma_1)$ a normalized Hecke eigenform of weight $2k$, 
$F_2\in S_{k+1}(\Gamma_2) =Maass(f)$ the Maass lift of $f$,  and in general
$F_{2n}\in S_{k+n}(\Gamma_{2n})$ the Ikeda lift of $f$
(we assume $k\equiv n\bmod 2$, $n\in \N$).
Next let
$g\in S_{k+n+r}(\Gamma_{r})$  be an arbitrary   Siegel cusp eigenform of genus $r$ and weight $k+n+r$,  with $n,r\ge 1$.Then according to Ikeda-Miyawaki there exists a Siegel eigenform
${\cal F}_{f,g}\in S_{k+n+r}(\Gamma_{2n+r})$
such that 
$$
L(s,{\cal F}_{f,g}, St) = L(s,g,St)
\prod_{j=1}^{2n}L(s+k+n-j,f) 
$$
(under a non-vanishing condition).

\subsection*{{$p$-adic versions of Ikeda-Miyawaki constructions}}
Now consider a $p$-adic family 
$$
k{}\mapsto f_k=\sum_{n=1}^\infty a_{n}(k)q^n \in {\overline{\mathbb{Q}}}[\![q]\!] \subset \mathbb{C}_p[\![q]\!]
$$ 
Fourier coefficients $a_n(k)$ of $f_k$ and one of  Satake $p$-parameters  $\alpha(k):=\alpha_p^{(1)}(k)$ are given by the 
certain analytical $p$-adic functions $k{}\mapsto a_n(k)$ for $(n,p)=1$.  

Then the Fourier expansions of the modular forms $F=F_k$ and ${\cal F}_{f,g}={\cal F}_{f_k,g}$
can be explicitly evaluated giving examples of $p$-adic Siegel modular forms of positive slope.

Notice that the Satake parameters of $F$ are 
$
\beta_0, \beta_1, \dots, \beta_{2{n}}$, where 
$$
\beta_0=p^{{n}k-{n}({n}+1)/2},
\beta_i=\alpha (k) p^{i-k}, \beta_{{n}+i}=\alpha(k)^{-1} p^{k+i-1} (i=1, \dots, {n}).
$$ 

\subsection*{Acknowledgement}

It is a great pleasure for me to thank  
Siegfried Boecherer, 
Stephen Gelbart,
Solomon Friedberg and Wadim 
Zudilin 
 for valuable  discussions and observations. 

\

My special thanks go to 
Yuri 
Nesterenko for the invitation to 
the International conference 
"Diophantine and Analytical Problems in Number Theory" 
for  100th birthday of  A.O. Gelfond in Moscow University, 
and for
soliciting a paper for the proceedings.


\bibliographystyle{plain}

\end{document}